\documentclass{article}
\usepackage[a4paper, margin=1in]{geometry} 
\usepackage{natbib}
\usepackage{amsmath, amssymb} 
\usepackage{wrapfig}
\usepackage{amsmath}
\usepackage{hyperref}
\usepackage{booktabs}
\usepackage{rotating}
\usepackage{multirow}
\usepackage{array}
\usepackage{graphicx}
\usepackage{caption}
\usepackage{siunitx} 
\usepackage{titlesec}
\usepackage{subcaption} 
\usepackage{empheq}
\usepackage{times}
\usepackage{geometry}

\title{Solving Newell-Whitehead-Segel and Allen-Cahn Equations Employing Physics-Informed Neural Networks: A Comparative Analysis with Spline Methods}
\author{
    Ali Haider Shah$^{1}$\thanks{E-mail: ah894@leicester.ac.uk (Corresponding author)} \and
    Naveed R. Butt$^{2}$\thanks{E-mail: naveedrbutt@gmail.com }\and
    Asif Ahmad$^{3}$\thanks{E-mail: asifahmad7007@gmail.com } \and
    Muhammad Omer Bin Saeed$^{4}$\thanks{E-mail: omer.saeed@giki.edu.pk}
}
\date{}

\begin{document}

\maketitle

\begin{center}

\textsuperscript{1}{School of Computing and Mathematical Sciences, University of Leicester, UK}

\textsuperscript{2, 3, 4}{Faculty of Engineering Sciences, Ghulam Ishaq Khan Institute of Engineering Sciences and Technology, Topi, Pakistan.}

\end{center}

\begin{abstract}
This study focuses on the solution of partial differential equations (PDEs) by using physics-informed neural networks (PINNs). The Newell-Whitehead-Segel (NWS) equation and the Allen-Cahn equation belong to fundamental PDEs used mostly in various scientific disciplines. Different methods, including analytical and numerical approaches, have been proposed for solving these equations alongside the recently introduced PINN method. This study provides a detailed and comprehensive comparison between the developed PINN method and the state-of-the-art spline numerical solution for the NWS and Allen-Cahn equation. Furthermore, the computational time of the trained PINN models is evaluated to determine their computational efficiency. The findings show that PINN is significantly better than spline methods in solving both problems.
\end{abstract}

\section{Introduction}
PDEs have the potential to explain more complicated processes arising in physics, chemistry, and biology. Therefore the evaluation of solutions to PDEs is a valuable and effective way of studying the dynamic behavior or the functions involved in such equations. The NWS equations describe the behavior of Rayleigh-Bénard convection close to the bifurcation point of binary fluid mixtures.\cite{rosu2005supersymmetric}. 
Rayleigh–Bénard convection is a specific type of natural convection present in a flat layer of fluid being heated from below, which causes the fluid to form a set of closed circulation cells named the cells of Bénard. There is always the creation of convective motion when the heating is enough high: hot fluid is located at the top and cold one is at the bottom. As Rayleigh-Bénard convection is relatively easy to study analytically and experimentally, it is investigated. One prominent illustration of self-organizing nonlinear systems is provided by convection patterns
\cite{getling1998rayleigh}. Two similar structures observed in Rayleigh-Bénard convection are the roll pattern or we can say stripe pattern and the hexagonal pattern  \cite {golovinb2007general}.
\\
The Allen-Cahn equation emerges in numerous scientific disciplines, including biomathematics, quantum mechanics, and plasma physics \cite{haq2024approximate}. The Allen–Cahn initially introduced the  Allan Cahn model 
in 1979. It is famous that plasma medium and fluid dynamics phenomena are successfully modeled by kink-shaped tanh solutions or bell-shaped sech solutions \cite{inan2020analytical}.

Since it is conventional to use approaches linked to the intrinsic properties of PDEs for constructing exact solutions,Therefore Sophisticated numerical techniques are used to find the numerical solution of PDEs, such as legendre wavelet-based approximation method \cite{hariharan2014efficient}, convergent finite difference method \cite{bates2009numerical},  Haar wavelet
 method \cite{hariharan2009haar}, invariant energy quadratization (IEQ)
method \cite{yang2020convergence}, non polynomial b-spline collocation method  \cite{haq2024approximate} and Reproducing kernel function \cite{xie2008numerical}. Specifically Allen Cahn equation of different orders and degrees has been successfully addressed by using aforementioned numerical methods.\\
 In the recent years, Physics Informed Neural Networks (PINN) has earned significant attention through the ground break work of Raissi. In his work he utilized neural networks to find the solution of PDE's. This novel approach, which falls under the deep learning domain, has attracted a lot of interest due to its capacity to facilitate data-driven PDE solution finding \cite{raissi2017physics}. So far, PINN has found extensive application in scientific computing, specifically in addressing both forward and inverse problem in non linear PDEs \cite{pu2021solving,zhang2023enforcing,blechschmidt2021three,wang2021data,kumar2024approximate}. Recently, PINN has been compared with other numerical methods in terms of accuracy and computational efficiency \cite{savovic2023comparative}.
\\
The approximate solution by PINN might not always be feasible and can show significant departures from the desired result. Observation suggests that a variety of reasons might be behind these differences. Neural networks incorporate many factors, and in order to obtain the best results, each parameter must be well-optimized. For this reason, we continuously modify and enhance our neural networks to get optimized solutions. This path of development has led to the modified versions of PINN. For instance, to address both forward and inverse PDE problems, for example, a gradient-enhanced PINN has been created, which incorporates gradient information from the PDEs residual into the loss function \cite{yu2022gradient}. Lin and Chen proposed
a two-stage PINN where the first stage is standard PINN and the second stage integrates the
conserved quantities into mean squared error loss to train neural networks \cite{lin2022two}. Zhu et al. built group equivariant neural networks that effectively simulated various periodic solutions of nonlinear dynamical lattices while adhering to spatiotemporal parity symmetries \cite{zhu2022neural}. To speed up training and convergence, Jagtag et al. presented local adaptive activation functions by adding a scaling parameter to each layer and neuron independently \cite{jagtap2020locally}. Some of useful techniques such as properly designed non-uniform training weighing technique is very useful for improving accuracy of PINN \cite{wang2021understanding}. To remove obstacles in more complex and realistic applications, more focused neural networks were developed, including the Bayesian PINN \cite{yang2021b}, the discrete PINN framework based on graph convolutional network and variational structure of PDE \cite{gao2022physics}, extended physics-informed neural networks (XPINNs) \cite{jagtap2020extended}, peripheral physics-informed neural network (PPINN) \cite{meng2020ppinn}, auxiliary physics-informed neural networks (A-PINN) \cite{yuan2022pinn}, fractional physics-informed neural networks (fPINNs) \cite{pang2019fpinns}, and variational PINN \cite{kharazmi2021hp}. To solve PDEs, neural networks require independent information. Leveraging some of the intrinsic properties of PDEs through Lie symmetry in PINN can be advantageous for enhancing the solution. Embedding symmetry loss information in total neural network loss has produced good results. \cite{li2023utilizing,zhang2023enforcing}.
\\
\\
The main idea  of this study is to compare PINNs with state of the art numerical methods (specifically Spline methods) and  enhance their performance. To achieve this goal two problems, NWS equation and  Allen Cahn equation  are examined. In this study, numerous parameters of PINN are tuned to achieve best results. The findings show that PINN performs remarkably better than other numerical techniques in terms of accuracy. Furthermore, computational efficiency is evaluated after training the model on 1,000 to 10,000 domain points, with intervals of 1,000. The computation time shows a very nice approximate linear behavior, indicating that the order of complexity is linear, i.e, O(n).

\section{Methodology}

The general form of a nonlinear parabolic PDE is
\begin{equation}
u_t = m u_{xx} + nu + ou^p + \eta(x,t,u,u_x)
, \quad \quad (x,t) \in X \times [0, T]  \label{eq:1}
\end{equation}
\begin{equation}
u(x,0) = f(x), \quad \quad   x \in X \subset (a,b)
\end{equation}

\begin{equation}
u(a, t) = g(t) \quad \text{and} \quad u(b, t) = h(t), \quad t \in [0, T].   
\end{equation}

These equations show the behavior of the solution \( u \) with respect to time \( t \) and the spatial coordinate \( x \), where \(m\), \(n\), and \(o\) are coefficients, and \( a \), \( b \), and \( T \) are constants. The term \(\eta(x,t,u,u_x)\) represents a general source term or nonlinear function that can depend on the spatial coordinate \(x\), the time \(t\), the solution \(u\), and its spatial derivative \(u_x\). This term can account for various effects such as external forces, reaction terms, or additional complex interactions. The function \( f(x) \) represents the initial value of \( u \), while \(g(t)\) and \(h(t)\) are functions of \( t \) at the boundary points.
 We examine a feed-forward Neural Network shown in Figure \ref{fig:PINN}. in order to comprehend the behaviour of a neural network. Consider a NN of \( L \) layers. This structure consists of one input layer, \( L - 1 \) hidden layers and one output layer. \(N_{l}\) are the number of neurons for each \( l \)-th layer. In this network, each layer receives an output \( Z_{l-1} \in \mathbb{R}^{N_{l-1}} \) from the previous layer, which goes through the affine transformation, followed by the application of non-linear activation function \( \sigma(\cdot) \).  defined mathematically in eqn (4).

\begin{equation}
Z_{l} = \sigma(W_{l} Z_{l-1} + b_{l}) 
\end{equation}

Here, \( W_{l} \in \mathbb{R}^{N_{l} \times N_{l-1}} \) and \( b_{l} \in \mathbb{R}^{N_{l}} \) are the weights and biases of the \( l \)-th layer, respectively. For output layer, we use different  activation functions according to our required result. Thus, the final output of the NN is given by eqn (5):

\begin{equation}
 O_{\theta}(Z) = Z_{L}    
\end{equation}
 
where \( Z \) = \( Z_{0}\)  represents the input, \( \theta = \{ W_l, b_l \}_{l=1}^{L} \) are the trainable parameters of the NN and \( O \) signifies the output of the NN. Furthermore we use different optimizer to find best possible result according to nature of our problem and defined loss function.
In general loss function can be defined as:

\begin{equation}
\text{loss} = \sum_{i=1}^{n} (y_i - O_{\theta}(Z_i))^2    
\end{equation}

where \(y_i\) represents the actual value and \(O_{\theta}(Z_i)\) represents the predicted value by NN for the \(i\)-th sample, and \(n\) is the total number of sample points.This is generalized structure to represent neural network its architecture changes according to nature of problem. Eqn (6) represents the data driven model. In the next section 2.1, we explain the variant of neural network called PINN that is used for solution of PDEs. PINN is known as blend of  model driven model along with data-driven model.
\\
\subsection{PDE's solution via PINN}
\begin{figure}[htbp] 
    \centering
    \includegraphics[width=0.9\textwidth]{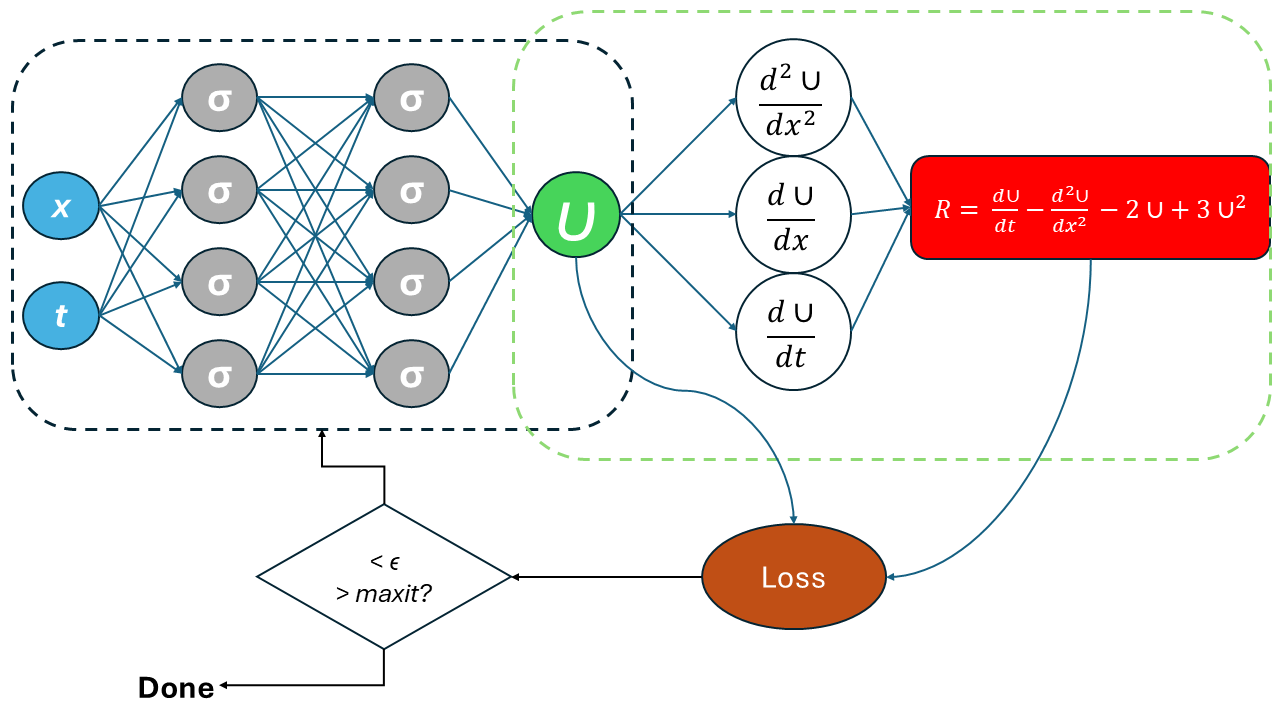} 
    \caption{Physics Informed Neural Network for NWS Equation}
    \label{fig:PINN} 
\end{figure}

The basic idea behind PINN for solving PDEs is to utilize a neural network $u_{\text{NN}}(t, x, \Theta)$ parameterized by $\Theta$ to approximate the exact solutions $u(t, x)$. In this technique, the neural network architecture is designed to encode the fundamental physics of the system, enabling it to learn the solution behavior directly from the data as shown in Figure\eqref{fig:PINN}.
\\
\\
Let's denote the data for the initial condition and boundary condition as \( \{(t_i^0, x_i^0, u_i^0)\}_{i=1}^{N_0} \) and \( \{(t_i^b, x_i^b, u_i^b)\}_{i=1}^{N_b} \), respectively, where \( N_0 \) and \( N_b \) represent the respective numbers of data points. Additionally, let \( \{(t_j^c, x_j^c)\}_{j=1}^{N_c} \) represent the collocation points used for evaluating the residual of the PDE, where \( N_c \) is the number of collocation points.
The loss term for the initial condition is denoted as \( \mathcal{L}_{\text{init}}(\Theta) \), for the boundary condition as \( \mathcal{L}_{\text{bound}}(\Theta) \), and for the residual term as \( \mathcal{L}_{\text{res}}(\Theta) \). These terms are defined as:

\begin{equation}
\mathcal{L}_{\text{init}}(\Theta) = \frac{1}{N_0} \sum_{i=1}^{N_0} \left( u_{\text{NN}}(t_i^0, x_i^0, \Theta) - u_i^0 \right)^2
\end{equation}

\begin{equation}
\mathcal{L}_{\text{bound}}(\Theta) = \frac{1}{N_b} \sum_{i=1}^{N_b} \left( u_{\text{NN}}(t_i^b, x_i^b, \Theta) - u_i^b \right)^2  
\end{equation}

\begin{equation}  \mathcal{L}_{\text{res}}(\Theta) = \frac{1}{N_c} \sum_{j=1}^{N_c} \left( \mathcal{R}(t_j^c, x_j^c, \Theta) \right)^2 
\end{equation}

where $u_{\text{NN}}(t_i^0, x_i^0, \Theta)$ ,$u_{\text{NN}}(t_i^b, x_i^b, \Theta)$ represents the predicted solution of PDE at initial and boundary points by NN. The residual of the PDE at collocation points denoted as
$\mathcal{R}(t_j^c, x_j^c, \Theta)$.
Residual term for eqn (1) defined as:
\begin{equation}
\mathcal{R} = u_t - mu_{xx} - nu - ou^p - \eta(x,t,u,u_x)
\end{equation}

\begin{figure}[ht]
    \centering
    \begin{subfigure}[b]{0.45\textwidth}
        \centering
        \includegraphics[width=\textwidth]{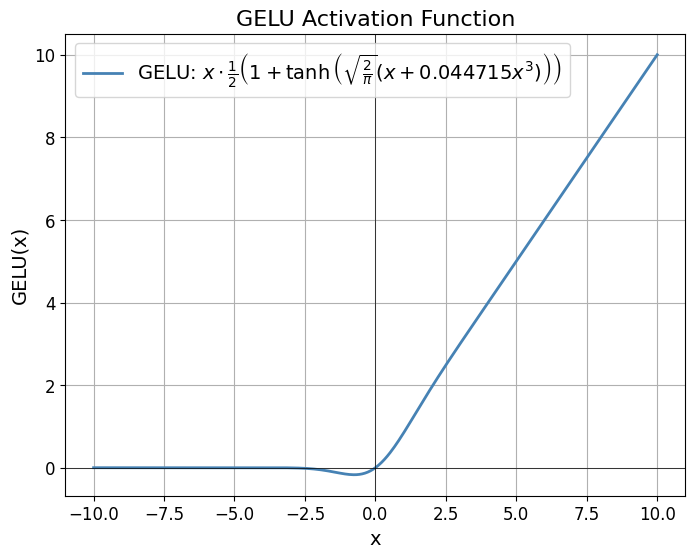}
        
        \label{fig:figure1}
    \end{subfigure}
    \hfill
    \begin{subfigure}[b]{0.45\textwidth}
        \centering
        \includegraphics[width=\textwidth]{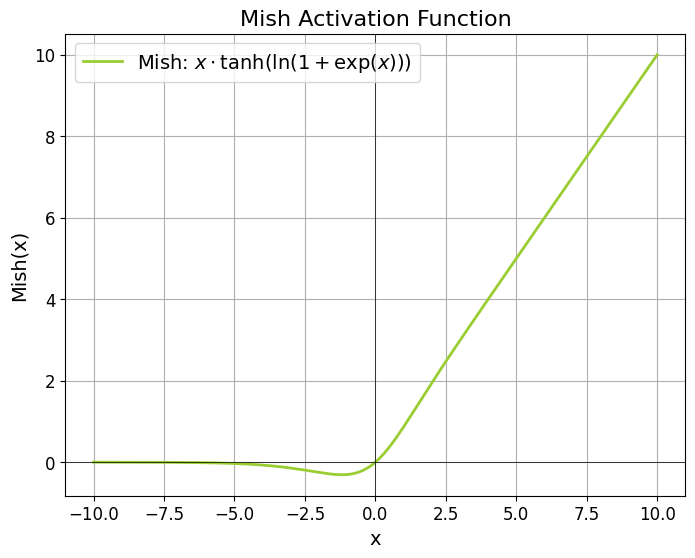}
        
        \label{fig:figure2}
    \end{subfigure}
    \vfill
    \begin{subfigure}[b]{0.45\textwidth}
        \centering
        \includegraphics[width=\textwidth]{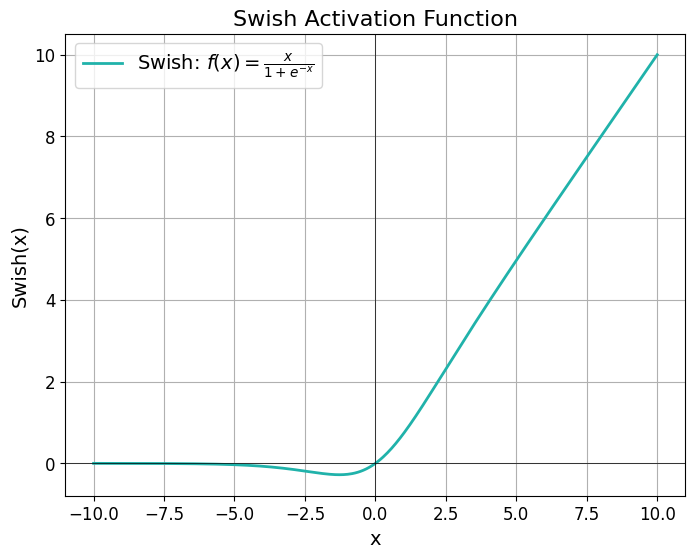}
        
        \label{fig:figure3}
    \end{subfigure}
    \hfill
    \begin{subfigure}[b]{0.45\textwidth}
        \centering
        \includegraphics[width=\textwidth]{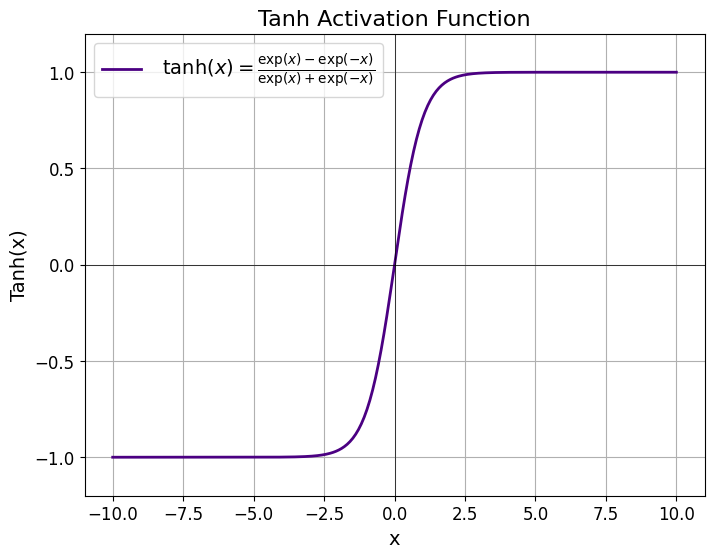}
        
        \label{fig:figure4}
    \end{subfigure}
    \caption{Activation functions}
    \label{fig:main}
\end{figure}

Finally  total loss function combining all terms in eqn(7), eqn(8) and eqn(9) can be written as:

\begin{equation}
 \mathcal{L}(\Theta) = \alpha \mathcal{L}_{\text{init}}(\Theta) + \beta \mathcal{L}_{\text{bound}}(\Theta) + \gamma \mathcal{L}_{\text{res}}(\Theta)    
\end{equation}

where $\alpha$, $\beta$ and $\gamma$ are weights, usually they are considered as $\alpha$ = $\beta$ = $\gamma$ = 1.
The PINN technique in this situation often makes use of the Adam \cite{kingma2014adam} or L-BFGS \cite{liu1989limited} algorithms to update the parameters $\Theta$ and, as a result, minimize the mean square error loss function. We can put constraints that contain Independent properties and information of  our  system to make neural network more efficient.
Furthermore, activation functions play an key role in the PINN methodology. A recent comparative study \cite{lima2022multilayer} on solving PDE models using PINNs demonstrated that various activation functions were utilized. The results showed that the hyperbolic tangent (tanh) and Gaussian Error Linear Unit (GELU) activation functions performed the best against all activation functions as shown in Fig \ref{fig:main}. In this research, GELU was found to work exceptionally well against all other activation functions, which is why GELU is used for both problems in this study.

\section{Numerical Experiment}

\subsection{Problem 1} 
\subsubsection*{Newell-Whitehead-Segel equation}

if \(n = 2\), \(o = -3\), \(p = 2\), \(m = 1\) in \eqref{eq:1} and \(\eta(x,t,u,u_x) = 0\) then the NWS equation with initial and boundary conditions is written  as \cite{zahra2014cubic}:

\begin{equation}
u_t - u_{xx} - 2u + 3u^2 = 0,
\end{equation}

\begin{equation}
u(x,0) = \lambda, \quad u(0,t) = \frac{-2\lambda e^{2t}}{-2 + 3\lambda(1 - e^{2t})},
\end{equation}

\begin{equation}
u(1,t) = \frac{-2\lambda e^{2t}}{-2 + 3\lambda(1 - e^{2t})}
\end{equation}
The general solution is

\begin{equation}
 u(x,t) = \frac{-2\lambda e^{2t}}{-2 + 3\lambda(1 - e^{2t})}
\end{equation}

\begin{table}[htbp]
  \centering
  \caption{Fixed constructional Components of the PINN for problem 1}
    \begin{tabular}{ll}
    \toprule
    \textbf{Component} & \textbf{Value} \\
    \midrule
    Activation function & GELU \\
    Number of layers & 8 \\
    Number of neurons per layer & 20 \\
    Initial condition points, $N_0$ & 250 \\
    Boundary condition points, $N_b$ & 250\\
    Collocation points, $N_r$ & 10,000 \\
    Learning rate schedule & Piece wise Decay \\
    & $\left[
    \begin{array}{c}
    1000 \\ 10^{-2}
    \end{array} \right] 
    \left[ \begin{array}{c}
    3000 \\ 10^{-3}
    \end{array} \right]
    \left[ \begin{array}{c}
    5000 \\ 5 \times 10^{-4}
    \end{array} \right]$ \\
    Optimizer & Adam \\
    Training (iterations) & 20,000 \\
    \bottomrule
  \end{tabular}
  \label{tab:fixed_componentsaa}
\end{table}

\begin{table}[htbp]
\centering
\caption{Absolute error for $0 \leq x \leq 1$, $0 \leq t \leq 1$ and $\lambda = 0.1$ \cite{zahra2014cubic}}
\begin{tabular}{|c|c|c|c|c|c|c|}
\hline
$x/t$ & Method & 0.2 & 0.4 & 0.6 & 0.8 & 1.0 \\ \hline
0.2 & UCBS & $8.323 \times 10^{-4}$ & $1.011 \times 10^{-3}$ & $8.581 \times 10^{-4}$ & $4.745 \times 10^{-4}$ & $1.991 \times 10^{-5}$ \\
\cline{2-7}
& TCBS & $8.295 \times 10^{-4}$ & $1.007 \times 10^{-3}$ & $8.514 \times 10^{-4}$ & $4.651 \times 10^{-4}$ & $3.239 \times 10^{-5}$ \\
\cline{2-7}
& ECBS & $6.068 \times 10^{-4}$ & $6.800 \times 10^{-4}$ & $5.476 \times 10^{-4}$ & $2.746 \times 10^{-4}$ & $6.339 \times 10^{-5}$ \\
\cline{2-7}
& \textbf{PINN}  & $6.270 \times 10^{-6}$ & $1.012 \times 10^{-5}$ & $3.110 \times 10^{-6}$ & $9.390 \times 10^{-6}$ & $3.960 \times 10^{-6}$ \\
\hline
0.4 & UCBS & $1.226 \times 10^{-3}$ & $1.520 \times 10^{-3}$ & $1.302 \times 10^{-3}$ & $7.334 \times 10^{-4}$ & $4.932 \times 10^{-6}$ \\ 
\cline{2-7}
& TCBS & $1.222 \times 10^{-3}$ & $1.514 \times 10^{-3}$ & $1.292 \times 10^{-3}$ & $7.194 \times 10^{-4}$ & $2.342 \times 10^{-5}$ \\
\cline{2-7}
& ECBS & $9.013 \times 10^{-4}$ & $1.023 \times 10^{-3}$ & $8.289 \times 10^{-4}$ & $4.221 \times 10^{-4}$ & $8.366 \times 10^{-5}$ \\
\cline{2-7}
& \textbf{PINN}  & $6.000 \times 10^{-8}$ & $1.143 \times 10^{-5}$ & $1.280 \times 10^{-6}$ & $1.404 \times 10^{-5}$ & $4.500 \times 10^{-6}$ \\
\hline
0.6 & UCBS & $1.226 \times 10^{-3}$ & $1.520 \times 10^{-3}$ & $1.302 \times 10^{-3}$ & $7.334 \times 10^{-4}$ & $4.932 \times 10^{-6}$ \\
\cline{2-7}
& TCBS & $1.222 \times 10^{-3}$ & $1.514 \times 10^{-3}$ & $1.292 \times 10^{-3}$ & $7.194 \times 10^{-4}$ & $2.342 \times 10^{-5}$ \\
\cline{2-7}
& ECBS & $9.013 \times 10^{-4}$ & $1.023 \times 10^{-3}$ & $8.289 \times 10^{-4}$ & $4.221 \times 10^{-4}$ & $8.366 \times 10^{-5}$ \\
\cline{2-7}
& \textbf{PINN} & $9.300 \times 10^{-6}$ & $1.232 \times 10^{-5}$ & $6.090 \times 10^{-6}$ & $1.588 \times 10^{-5}$ & $1.040 \times 10^{-5}$ \\
\hline
0.8 & UCBS & $8.323 \times 10^{-4}$ & $1.011 \times 10^{-3}$ & $8.581 \times 10^{-4}$ & $4.745 \times 10^{-4}$ & $1.991 \times 10^{-5}$ \\
\cline{2-7}
& TCBS & $8.295 \times 10^{-4}$ & $1.007 \times 10^{-3}$ & $8.514 \times 10^{-4}$ & $4.651 \times 10^{-4}$ & $3.239 \times 10^{-5}$ \\
\cline{2-7}
& ECBS & $6.068 \times 10^{-4}$ & $6.800 \times 10^{-4}$ & $5.476 \times 10^{-4}$ & $2.746 \times 10^{-4}$ & $6.339 \times 10^{-5}$ \\
\cline{2-7}
& \textbf{PINN} & $2.298 \times 10^{-5}$ & $1.061 \times 10^{-5}$ & $1.949 \times 10^{-5}$ & $1.231 \times 10^{-5}$ & $1.004 \times 10^{-5}$ \\
\hline
\end{tabular}
\label{point-wise absolue error of p2}
\end{table}

Table \ref{tab:fixed_componentsaa} represents the constructional components of the PINN used for the NWS equation. This table shows that the PINN model is trained for 20,000 iterations with the Gelu activation function, where the learning rate is piecewise decayed. Specifically, for the first 1,000 iterations, the learning rate is $10^{-2}$; from 1,000 to 3,000 iterations, it is $10^{-3}$; and from 3,000 iterations onwards, it is $5 \times 10^{-4}$. The domain of equation is  $0\leq x\leq 1$ and $0\leq t\leq 1$ and $\lambda = 0.1$. Furthermore Table \ref{point-wise absolue error of p2} represents the absolute error comparison of uniform cubic B-spline (UCBS), extended cubic uniform B-spline (ECBS),trigonometric cubic B-spline (TCBS) and PINN method. 
\begin{equation}
    \text{Absolute Error} = |\text{Exact} - \text{Approximated}|
\end{equation}

\begin{figure}[htbp]
    \centering
    \begin{subfigure}[b]{0.45\textwidth}
        \centering
        \includegraphics[width=\textwidth]{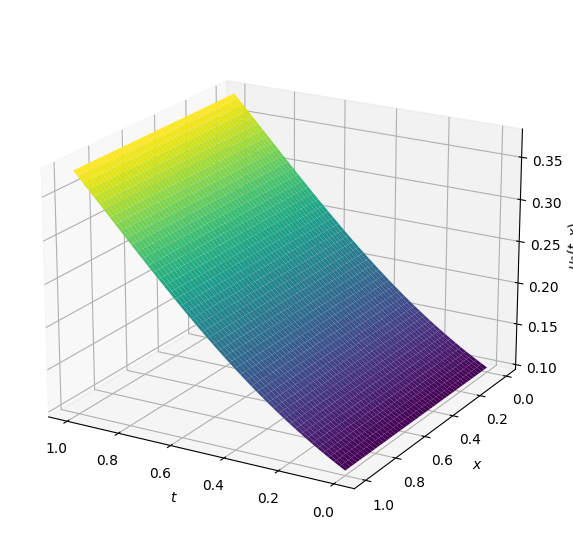} 
        \caption{Exact method}
        \label{exact 3d_1}
    \end{subfigure}
    \hfill
    \begin{subfigure}[b]{0.45\textwidth}
        \centering
        \includegraphics[width=\textwidth]{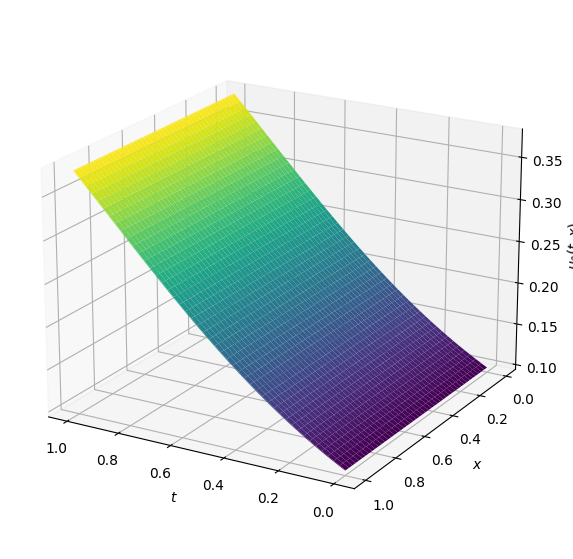} 
        \caption{PINN method}
        \label{PINN solution.png}
    \end{subfigure}
    \caption{ 3D surfaces of exact solution (left) and PINN solution(right) for Problem 1 where $\lambda = 0.1$}
    \label{3D_1}
    \end{figure}

\begin{figure}[htbp]
    \centering
    \begin{subfigure}[b]{0.45\textwidth}
        \centering
        \includegraphics[width=\textwidth]{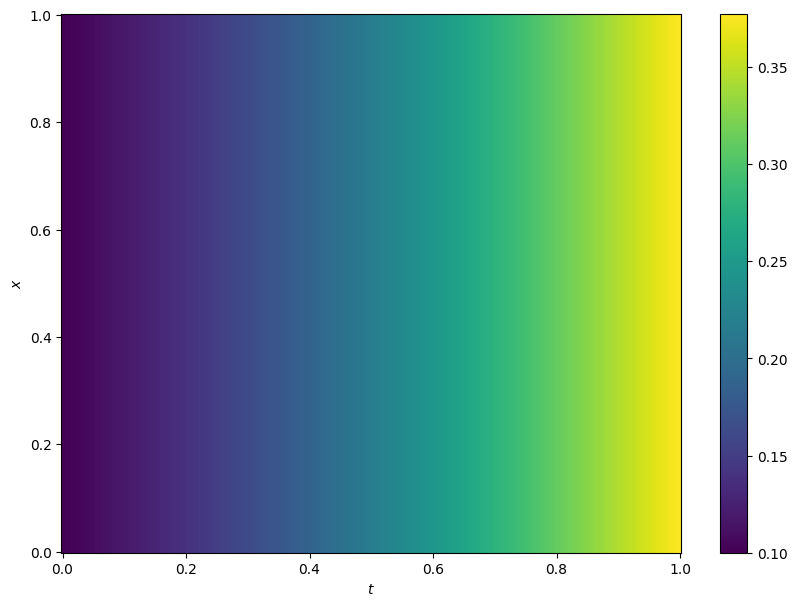} 
        \caption{Exact method}
        \label{heat map exact}
    \end{subfigure}
    \hfill
    \begin{subfigure}[b]{0.45\textwidth}
        \centering
        \includegraphics[width=\textwidth]{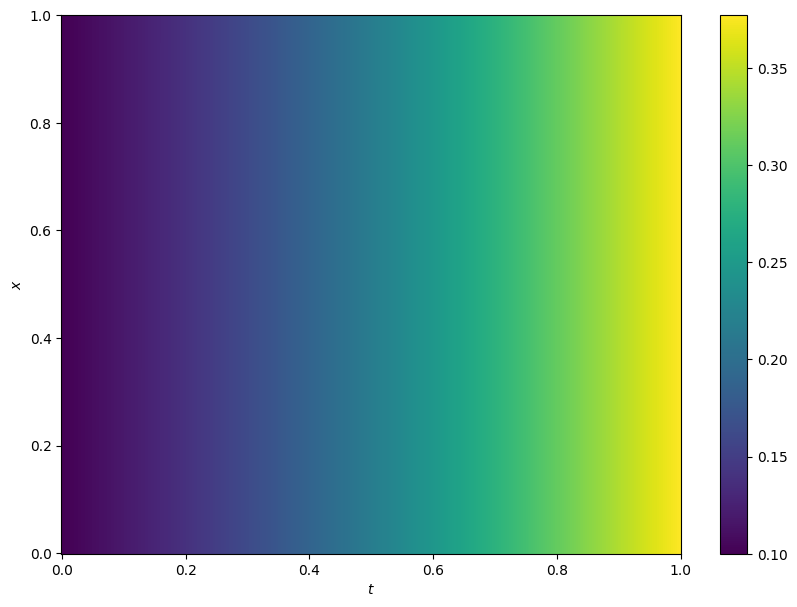} 
        \caption{PINN method}
        \label{heat map PINN}
    \end{subfigure}
    \caption{Heat maps plots representing the surface solutions for Problem 1 using both the exact method (left) and the PINN method (right)}
    \label{Heat map}
    \end{figure}

\begin{figure}[htbp]
    \centering
    \begin{subfigure}[b]{0.45\textwidth}
        \centering
        \includegraphics[width=\textwidth]{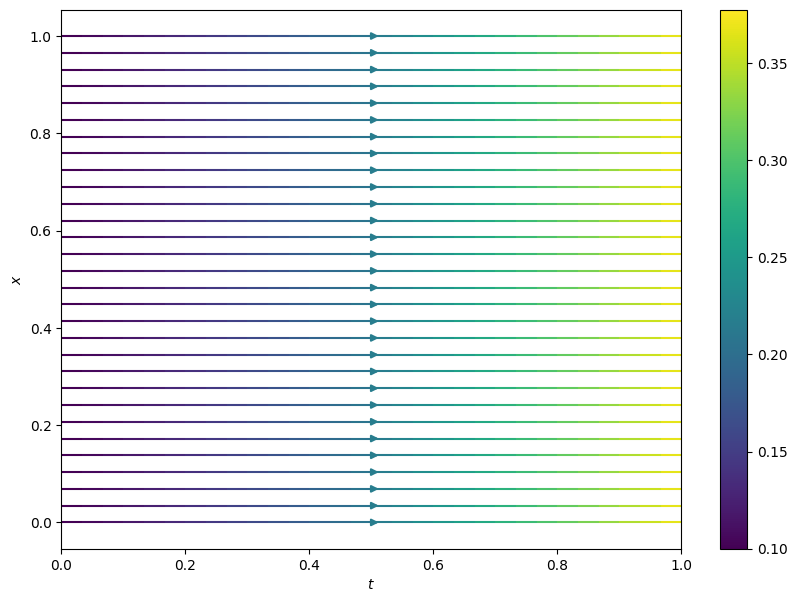} 
        \caption{Exact method}
        \label{streamline map exact}
    \end{subfigure}
    \hfill
    \begin{subfigure}[b]{0.45\textwidth}
        \centering
        \includegraphics[width=\textwidth]{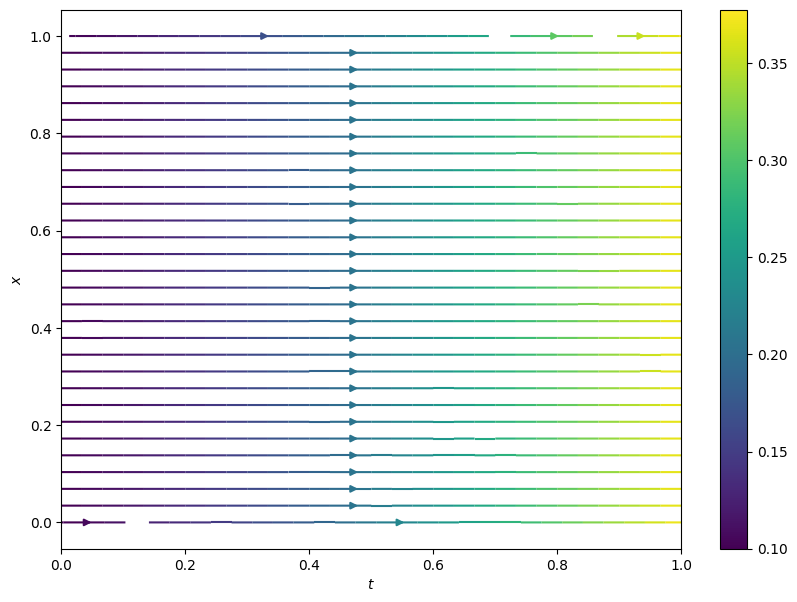} 
        \caption{PINN method}
        \label{u_pinn.png}
    \end{subfigure}
    \caption{Gradient lines plots representing the flow of gradient vectors at solution surface. where gradient lines of exact solution is on left side and the PINN method is on right side}
    \label{streamline map}
    \end{figure}

\begin{figure}[htbp]
    \centering
    \begin{subfigure}[b]{0.45\textwidth}
        \centering
        \includegraphics[width=\textwidth]{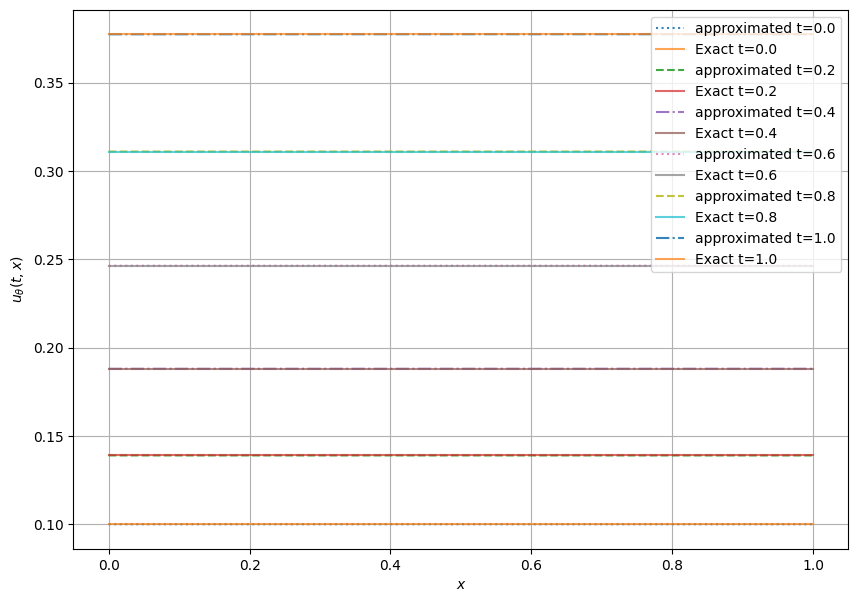} 
        \caption{Cross-section curves (2D)}
        \label{Crosssection curves}
    \end{subfigure}
    \hfill
    \begin{subfigure}[b]{0.45\textwidth}
        \centering
        \includegraphics[width=\textwidth]{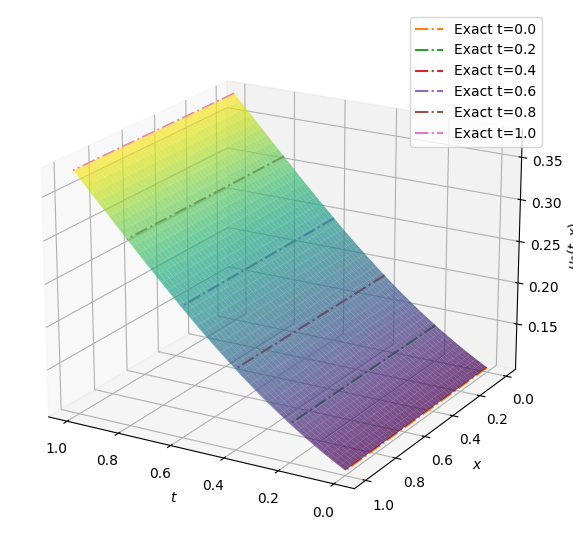} 
        \caption{Cross-section curves (3D)}
        \label{L2 norm loss}
    \end{subfigure}
    \caption{Comparison of cross-section curves (problem 1) between the exact method and the PINN method for various values of $t$ = 0, 0.2, 0.6, 0.8, 1. Dotted lines represent the exact solution cross-section curves, which fit well on the 3D solution surface generated by the PINN method (right). The 2D visualization of it, is shown on the left.}

    \label{loss}
    \end{figure}

\begin{figure}[htbp]
    \centering
    \begin{subfigure}[b]{0.45\textwidth}
        \centering
        \includegraphics[width=\textwidth]{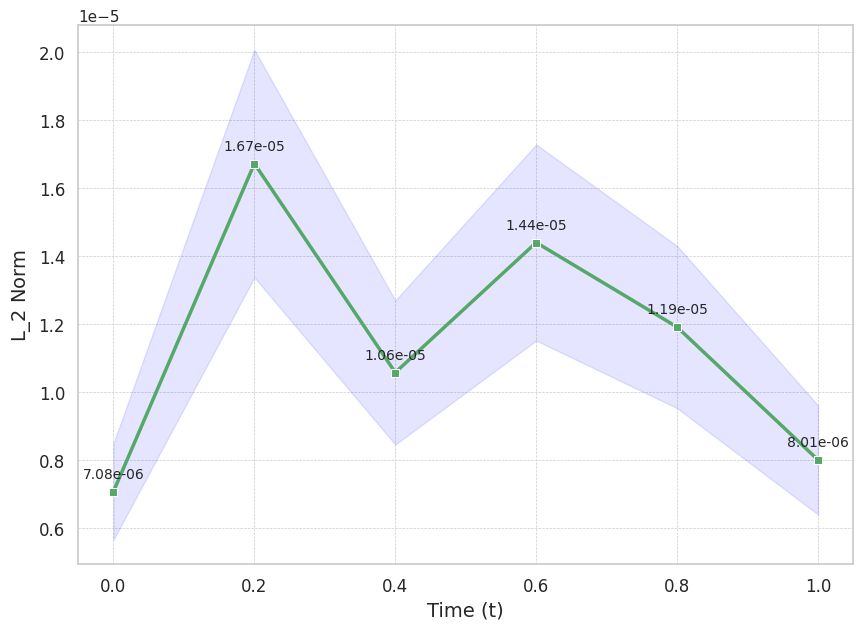} 
        \caption{\(L_2\) norm}
        \label{L_2 norm}
    \end{subfigure}
    \hfill
    \begin{subfigure}[b]{0.45\textwidth}
        \centering
        \includegraphics[width=\textwidth]{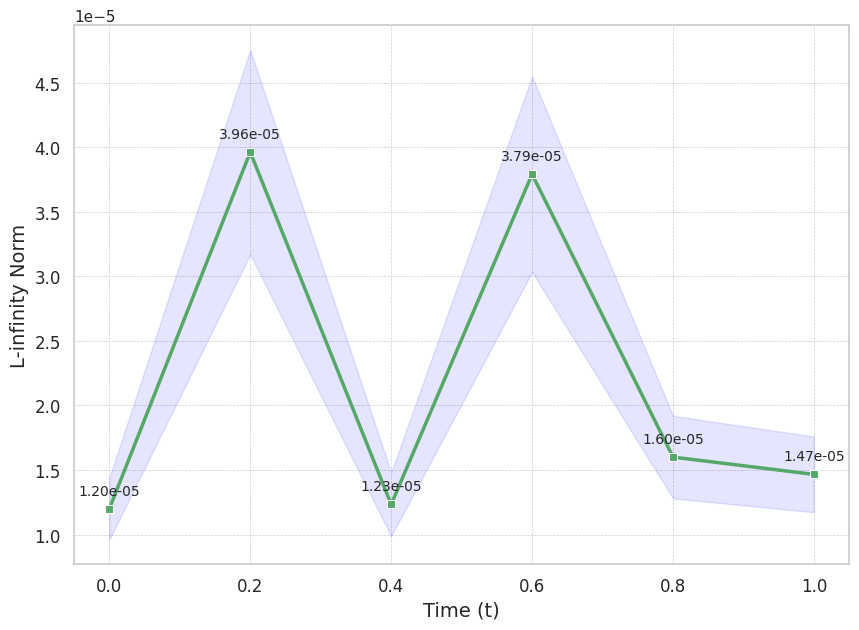} 
        \caption{$L_{\infty} norm$}
        \label{L_infty norm}
    \end{subfigure}
    \caption{ It represents the \(L_2\) norm and $L_{\infty}$ norm for problem 1 at t = 0, 0.2, 0.4, 0.6, 0.8, 1  by using PINN, where $ h = 0.004$}
    \label{L2 and L_infinity}
    \end{figure}

It is clear from Table \ref{point-wise absolue error of p2} that the PINN significantly outperforms all variants of spline numerical methods in terms of accuracy at the given points in the literature. Figure \ref{3D_1} represents the exact solution $u(x,t)$ and the approximate solution $u(x,t,\theta)$ by the PINN method, while Figure \ref{Heat map} represents the 2D visualization of Figure \ref{3D_1} in the form of a heat map.Figure \ref{streamline map} represents the flow of gradient vectors for the exact solution and the approximate solution by the PINN method. These flows of gradient vectors are represented by lines called gradient lines. Figure \ref{loss} represents the cross-section curves between exact and approximated values by PINN at different values of $t$. These cross-section curves of PINN overlap with the lines of the exact solutions, demonstrating how well PINN fits the exact solution. Figure \ref{L2 and L_infinity} depicts the $L_2$ and $L_\infty$ norms at various values of $t$, illustrating that accuracy remains consistently high over time.
where,
\begin{equation}
    L_2 = \sqrt{\sum_{i=1}^{n} (Exact_i - Approximated_i)^2}
\end{equation}

\begin{equation}
    L_\infty = \max_{i} |Exact_i - Approximated_i|
\end{equation}

\begin{figure}[htbp]
    \centering
\includegraphics[width=\textwidth]{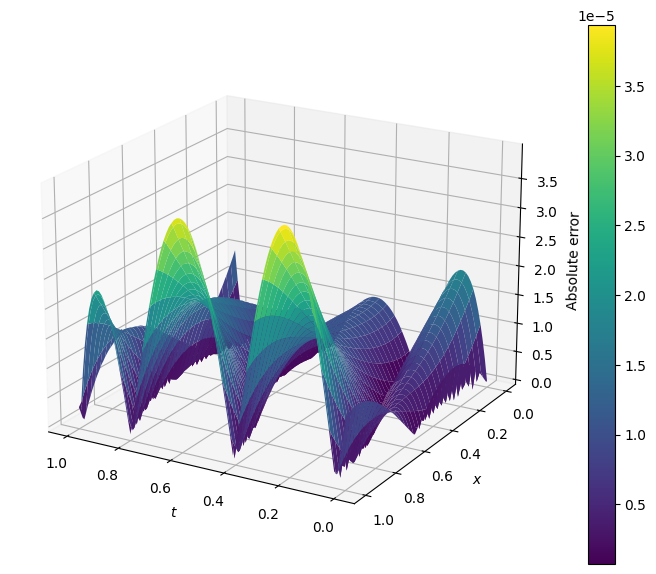} 
        \caption{Absolute error surface between exact and PINN method, where  $h = \delta t = 0.004$}
        \label{absolute_error_surface_h_delta_t_0.004}
    \end{figure}

\begin{figure}[htbp]
    \centering
    \begin{minipage}{0.45\textwidth}
        \centering
        \begin{tabular}{|c|c|}
            \hline
            \textbf{\shortstack{Number of \\ Points}} & \textbf{\shortstack{Computation Time \\ (seconds)}} \\
            \hline
            1000  & 10.812786  \\
            \hline
            2000  & 21.548918  \\
            \hline
            3000  & 31.382378  \\
            \hline
            4000  & 42.354686  \\
            \hline
            5000  & 54.906944  \\
            \hline
            6000  & 63.971302  \\
            \hline
            7000  & 74.595710  \\
            \hline
            8000  & 85.686710  \\
            \hline
            9000  & 97.562439  \\
            \hline
            10000 & 107.202053 \\
            \hline
        \end{tabular}
        \caption{Number of data points and corresponding computation times for problem 1.}
        \label{table:computation_time2}
    \end{minipage}
    \hfill
    \begin{minipage}{0.45\textwidth}
        \centering
        \includegraphics[width=\textwidth]{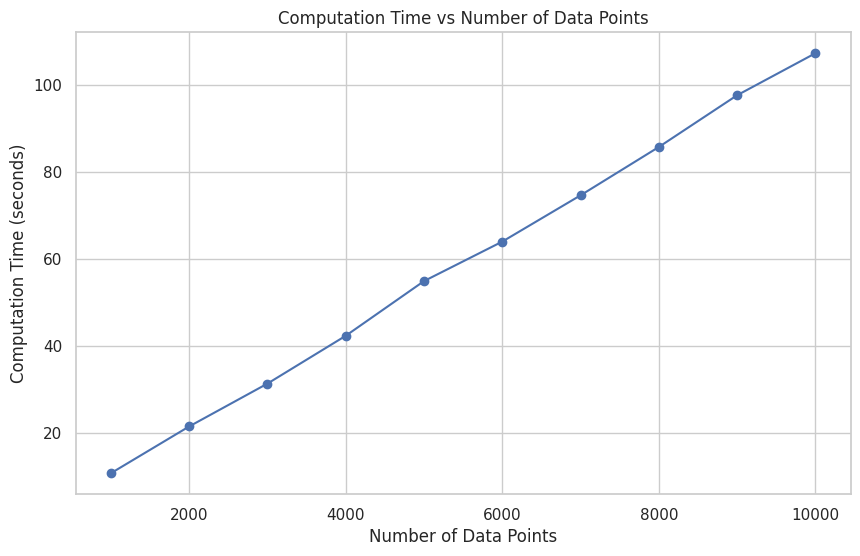}
        \caption{Plot of computation time vs. number of data points.}
        \label{fig:computation_time_plot2}
    \end{minipage}
\end{figure}

Figure \ref{absolute_error_surface_h_delta_t_0.004} displays the absolute error surface between the exact and PINN methods, where the spatial step size is represented by $h = 0.004$ and the time step size is denoted by $\delta t = 0.004$. Figure \ref{fig:computation_time_plot2} illustrates the computation time. After training the model, it evaluates the solution at random points within the domain, starting from 1000 and up to 10,000 with intervals of 1000. The results demonstrate that the model can compute the solution for 10,000 domain points in approximately 107 seconds, indicating computational efficiency. Furthermore, from Figure \ref{fig:computation_time_plot2}, it is evident that the computational complexity is linear, O(n), which signifies a robust model performance.

\newpage
\subsection{Problem 2}
\subsubsection*{Allen-Cahn Equation}

Here, we consider equation \eqref{eq:1} with parameters $m = 1, n = 1, o = -1, p = 3$ and \(\eta(x,t,u,u_x) = 0\) . The equation is defined in domain $0 < x < 1$ and $t > 0$ given by \cite{haq2024approximate},

\begin{equation}
u_t(x,t) - u_{xx}(x,t) - u(x,t) = -u^3,
\end{equation}

with the following conditions:

\begin{equation}
u(x,0) = -\frac{1}{2} + \frac{1}{2} \tanh(0.3536x),
\end{equation}

\begin{equation}
u(0,t) = -\frac{1}{2} + \frac{1}{2} \tanh(-0.75t),
\end{equation}

\begin{equation}
u(1,t) = -\frac{1}{2} + \frac{1}{2} \tanh(0.3536 - 0.75t).
\end{equation}

\begin{figure}[htbp]
    \centering
\includegraphics[width=\textwidth]{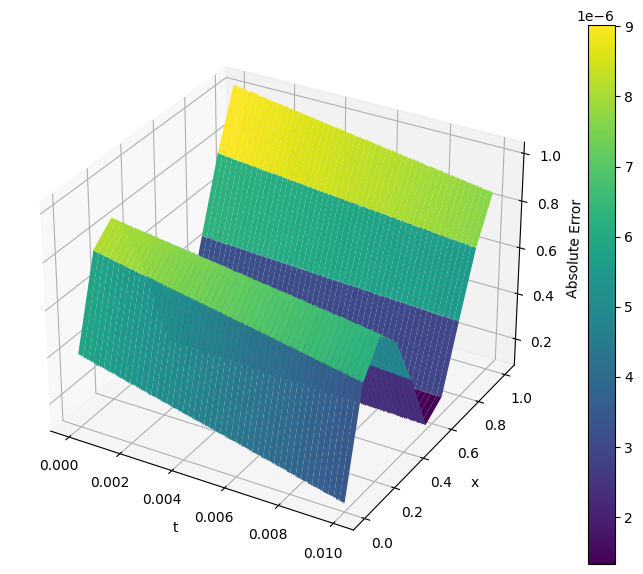} 
        \caption{Absolute error surface for problem 2, where $h = 0.1$ and $\delta t = 0.001$}
        \label{loss_function.png}
    \end{figure}

\begin{figure}[htbp]
    \centering
    \begin{minipage}{0.45\textwidth}
        \centering
        \begin{tabular}{|c|c|}
            \hline
            \textbf{\shortstack{Number of \\ Points}} & \textbf{\shortstack{Computation Time \\ (seconds)}} \\
            \hline
            1000  & 10.847024  \\
            \hline
            2000  & 22.224515  \\
            \hline
            3000  & 31.704301  \\
            \hline
            4000  & 42.410030  \\
            \hline
            5000  & 53.025968  \\
            \hline
            6000  & 64.038677   \\
            \hline
            7000  & 74.952875  \\
            \hline
            8000  & 83.142183  \\
            \hline
            9000  & 96.960017  \\
            \hline
            10000 & 104.414277 \\
            \hline
        \end{tabular}
        \caption{Number of data points and corresponding computation times for problem 2.}
        \label{table:computation_time for Allen cahn}
    \end{minipage}
    \hfill
    \begin{minipage}{0.45\textwidth}
        \centering
        \includegraphics[width=\textwidth]{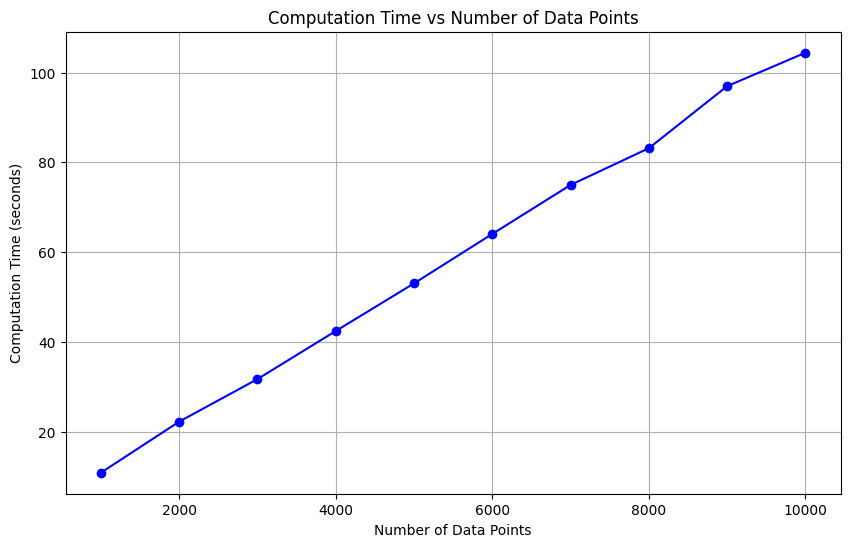 }
        \caption{Plot of computation time vs. number of data points.}
        \label{fig:computation_time for Allen }
    \end{minipage}
\end{figure}

The exact solution for the model equation is taken from:

\begin{equation}
u(x,t) = -\frac{1}{2} + \frac{1}{2} \tanh(0.3536x - 0.75t).
\end{equation}
The problem 2 has been addressed by the non-polynomial (N-P) spline technique and trigonometric cubic b-spline collocation (TCB-CM) method.
Table \ref{tab:PINN_absolute error} represents the point wise absolute error comparison, where , $\delta t$ and $h$
are step-size of time and spatial dimension \cite{haq2024approximate}. It is clear from  \ref{tab:PINN_absolute error} that PINN is far  more better than spline techniques. In this example, the PINN has better accuracy at all points given in the literature except at the boundary points. Figure \ref{loss_function.png} represents the absolute error surface where $0 \leq t \leq 0.01$ and $0 \leq x \leq 1$, for comparison with the cited paper. After training the model, the computation is also evaluated. The computation time for 10,000 points is approximately 104 seconds as shown in Figure \ref{fig:computation_time for Allen }.

\begin{table}[htbp]
  \centering
  \caption{Fixed Components of the PINN for problem 2 b}
    \begin{tabular}{ll}
    \toprule
    \textbf{Component} & \textbf{Value} \\
    \midrule
    Activation function & GELU \\
    Number of layers & 8 \\
    Number of neurons per layer & 40 \\
    Initial condition points, $N_0$ & 500 \\
    Boundary condition points, $N_b$ & 500\\
    Collocation points, $N_r$ & 10,000 \\
    Learning rate schedule & Piecewise Decay \\
    Optimizer & Adam  \\
    \bottomrule
    \end{tabular}%
  \label{tab:fixed_components of Allen Chan}%
\end{table}

\begin{table}[htbp]
\centering
\caption{Absolute error comparison of PINN method with present and [2], where $0 \leq t \leq 0.01$, $ 0 \leq x \leq 1 $, $\delta t = 0.001$ and $h = 0.1$}.
\begin{tabular}{|c|c|c|c|c|c|c|}
\hline
$\textbf{x/t}$ & 0.001 & 0.003 & 0.005 & 0.007 & 0.009 & 0.01 \\ \hline
\multicolumn{7}{|c|}{\textbf{N-P splines} \cite{haq2024approximate}} \\ \hline
0   & 0.00 & 0.00 & 0.00 & 0.00 & 0.00 & 0.00 \\ \hline
0.1 & 1.35$\times 10^{-4}$ & 3.79$\times 10^{-4}$ & 5.96$\times 10^{-4}$ & 7.92$\times 10^{-4}$ & 9.73$\times 10^{-4}$ & 1.06$\times 10^{-3}$ \\ \hline
0.2 & 9.12$\times 10^{-5}$ & 2.81$\times 10^{-4}$ & 4.74$\times 10^{-4}$ & 6.69$\times 10^{-4}$ & 8.63$\times 10^{-4}$ & 9.59$\times 10^{-4}$ \\ \hline
0.3 & 6.28$\times 10^{-5}$ & 1.91$\times 10^{-4}$ & 3.24$\times 10^{-4}$ & 4.61$\times 10^{-4}$ & 6.03$\times 10^{-4}$ & 6.75$\times 10^{-4}$ \\ \hline
0.4 & 3.64$\times 10^{-5}$ & 1.12$\times 10^{-4}$ & 1.91$\times 10^{-4}$ & 2.74$\times 10^{-4}$ & 3.61$\times 10^{-4}$ & 4.06$\times 10^{-4}$ \\ \hline
0.5 & 1.30$\times 10^{-5}$ & 4.14$\times 10^{-5}$ & 7.28$\times 10^{-5}$ & 1.07$\times 10^{-4}$ & 1.45$\times 10^{-4}$ & 1.65$\times 10^{-4}$ \\ \hline
0.6 & 7.69$\times 10^{-6}$ & 2.11$\times 10^{-5}$ & 3.18$\times 10^{-5}$ & 3.99$\times 10^{-5}$ & 4.52$\times 10^{-5}$ & 4.69$\times 10^{-5}$ \\ \hline
0.7 & 2.58$\times 10^{-5}$ & 7.57$\times 10^{-5}$ & 1.23$\times 10^{-4}$ & 1.69$\times 10^{-4}$ & 2.12$\times 10^{-4}$ & 2.33$\times 10^{-4}$ \\ \hline
0.8 & 4.12$\times 10^{-5}$ & 1.23$\times 10^{-4}$ & 2.04$\times 10^{-4}$ & 2.80$\times 10^{-4}$ & 3.53$\times 10^{-4}$ & 3.87$\times 10^{-4}$ \\ \hline
0.9 & 5.92$\times 10^{-5}$ & 1.64$\times 10^{-4}$ & 2.54$\times 10^{-4}$ & 3.32$\times 10^{-4}$ & 4.01$\times 10^{-4}$ & 4.32$\times 10^{-4}$ \\ \hline
1   & 0.00 & 0.00 & 0.00 & 0.00 & 0.00 & 0.00 \\ \hline
\multicolumn{7}{|c|}{\textbf{TCB-CM \cite{zahra2017trigonometric}}} \\ \hline
0   & 0.00 & 0.00 & 0.00 & 1.11$\times 10^{-16}$ & 1.11$\times 10^{-16}$ & 0.00 \\ \hline
0.1 & 2.45$\times 10^{-4}$ & 6.44$\times 10^{-4}$ & 9.69$\times 10^{-4}$ & 1.25$\times 10^{-3}$ & 1.49$\times 10^{-3}$ & 1.61$\times 10^{-3}$ \\ \hline
0.2 & 2.00$\times 10^{-4}$ & 6.13$\times 10^{-4}$ & 1.02$\times 10^{-3}$ & 1.40$\times 10^{-3}$ & 1.76$\times 10^{-3}$ & 1.94$\times 10^{-3}$ \\ \hline
0.3 & 1.80$\times 10^{-4}$ & 5.41$\times 10^{-4}$ & 9.08$\times 10^{-4}$ & 1.28$\times 10^{-3}$ & 1.64$\times 10^{-3}$ & 1.83$\times 10^{-3}$ \\ \hline
0.4 & 1.59$\times 10^{-4}$ & 4.81$\times 10^{-4}$ & 8.05$\times 10^{-4}$ & 1.13$\times 10^{-3}$ & 1.47$\times 10^{-3}$ & 1.63$\times 10^{-3}$ \\ \hline
0.5 & 1.41$\times 10^{-4}$ & 4.25$\times 10^{-4}$ & 7.12$\times 10^{-4}$ & 1.00$\times 10^{-3}$ & 1.30$\times 10^{-3}$ & 1.44$\times 10^{-3}$ \\ \hline
0.6 & 1.24$\times 10^{-4}$ & 3.75$\times 10^{-4}$ & 6.28$\times 10^{-4}$ & 8.84$\times 10^{-4}$ & 1.14$\times 10^{-3}$ & 1.27$\times 10^{-3}$ \\ \hline
0.7 & 1.09$\times 10^{-4}$ & 3.29$\times 10^{-4}$ & 5.52$\times 10^{-4}$ & 7.77$\times 10^{-4}$ & 1.00$\times 10^{-3}$ & 1.12$\times 10^{-3}$ \\ \hline
0.8 & 9.50$\times 10^{-5}$ & 2.90$\times 10^{-4}$ & 4.82$\times 10^{-4}$ & 6.71$\times 10^{-4}$ & 8.54$\times 10^{-4}$ & 9.43$\times 10^{-4}$ \\ \hline
0.9 & 8.89$\times 10^{-5}$ & 2.41$\times 10^{-4}$ & 3.73$\times 10^{-4}$ & 4.91$\times 10^{-4}$ & 6.01$\times 10^{-4}$ & 6.53$\times 10^{-4}$ \\ \hline
1   & 0.00 & 0.00 & 0.00 & 0.00 & 0.00 & 0.00 \\ \hline
\multicolumn{7}{|c|}{\textbf{ PINN}} \\ \hline
0   & 2.32$\times 10^{-6}$ & 2.56$\times 10^{-6}$ & 2.74$\times 10^{-6}$ & 2.98$\times 10^{-6}$ & 3.10$\times 10^{-6}$ & 3.22$\times 10^{-6}$ \\ \hline
0.1 & 1.85$\times 10^{-6}$ & 1.52$\times 10^{-6}$ & 1.34$\times 10^{-6}$ & 1.01$\times 10^{-6}$ & 9.54$\times 10^{-7}$ & 8.05$\times 10^{-7}$ \\ \hline
0.2 & 3.90$\times 10^{-6}$ & 3.70$\times 10^{-6}$ & 3.55$\times 10^{-6}$ & 3.31$\times 10^{-6}$ & 3.13$\times 10^{-6}$ & 3.04$\times 10^{-6}$ \\ \hline
0.3 & 4.47$\times 10^{-6}$ & 4.29$\times 10^{-6}$ & 4.08$\times 10^{-6}$ & 3.96$\times 10^{-6}$ & 3.84$\times 10^{-6}$ & 3.73$\times 10^{-6}$ \\ \hline
0.4 & 3.67$\times 10^{-6}$ & 3.55$\times 10^{-6}$ & 3.46$\times 10^{-6}$ & 3.34$\times 10^{-6}$ & 3.31$\times 10^{-6}$ & 3.16$\times 10^{-6}$ \\ \hline
0.5 & 1.97$\times 10^{-6}$ & 1.97$\times 10^{-6}$ & 1.85$\times 10^{-6}$ & 1.82$\times 10^{-6}$ & 1.79$\times 10^{-6}$ & 1.70$\times 10^{-6}$ \\ \hline
0.6 & 2.38$\times 10^{-7}$ & 2.68$\times 10^{-7}$ & 2.98$\times 10^{-7}$ & 3.58$\times 10^{-7}$ & 3.58$\times 10^{-7}$ & 3.28$\times 10^{-7}$ \\ \hline
0.7 & 2.32$\times 10^{-6}$ & 2.35$\times 10^{-6}$ & 2.41$\times 10^{-6}$ & 2.35$\times 10^{-6}$ & 2.38$\times 10^{-6}$ & 2.38$\times 10^{-6}$ \\ \hline
0.8 & 3.78$\times 10^{-6}$ & 3.78$\times 10^{-6}$ & 3.87$\times 10^{-6}$ & 3.87$\times 10^{-6}$ & 3.81$\times 10^{-6}$ & 3.87$\times 10^{-6}$ \\ \hline
0.9 & 4.38$\times 10^{-6}$ & 4.44$\times 10^{-6}$ & 4.44$\times 10^{-6}$ & 4.53$\times 10^{-6}$ & 4.44$\times 10^{-6}$ & 4.53$\times 10^{-6}$ \\ \hline
1   & 4.95$\times 10^{-6}$ & 4.98$\times 10^{-6}$ & 5.01$\times 10^{-6}$ & 5.07$\times 10^{-6}$ & 5.10$\times 10^{-6}$ & 5.10$\times 10^{-6}$ \\ \hline
\end{tabular}
\label{tab:PINN_absolute error}
\end{table}

\begin{figure}[htbp]
    \centering
    \begin{subfigure}[b]{0.45\textwidth}
        \centering
        \includegraphics[width=\textwidth]{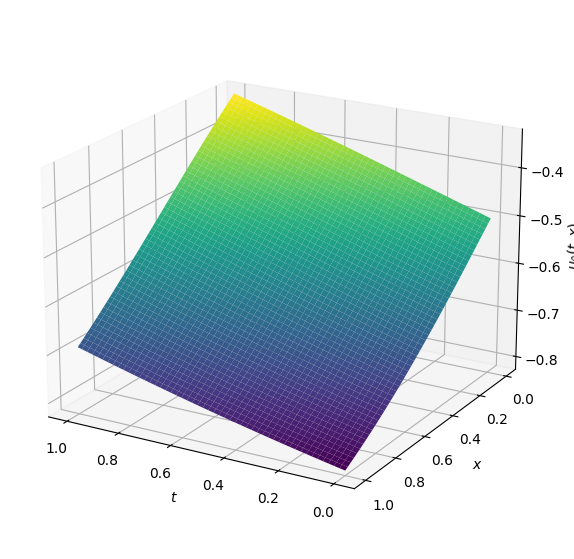} 
        \caption{Exact method}
        \label{ExactAllenchan.png}
    \end{subfigure}
    \hfill
    \begin{subfigure}[b]{0.45\textwidth}
        \centering
        \includegraphics[width=\textwidth]{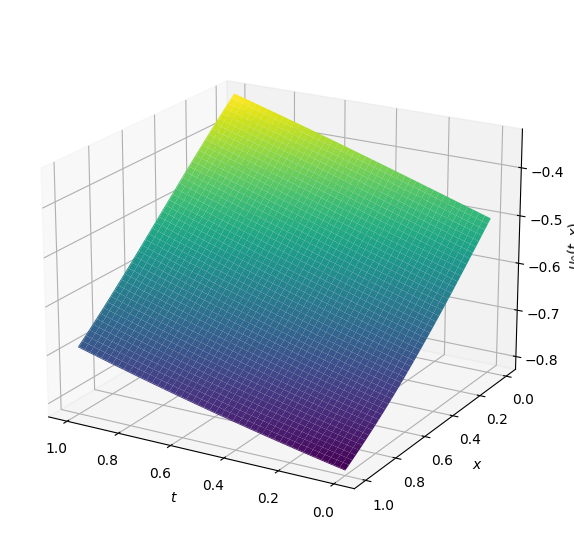} 
        \caption{PINN method}
        \label{ApproAllenchan.png}
    \end{subfigure}
    \caption{ 3D surfaces of exact solution (left) and PINN solution(right) for Problem 2}
    \label{3D_2}
    \end{figure}

\begin{figure}[htbp]
    \centering
    \begin{subfigure}[b]{0.45\textwidth}
        \centering
        \includegraphics[width=\textwidth]{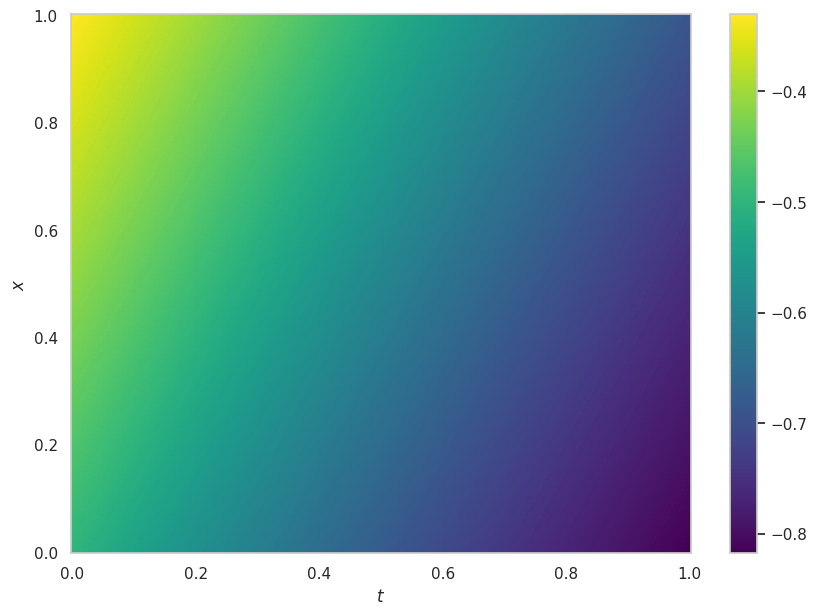} 
        \caption{Exact method}
        \label{Heatmap of exact.png}
    \end{subfigure}
    \hfill
    \begin{subfigure}[b]{0.45\textwidth}
        \centering
        \includegraphics[width=\textwidth]{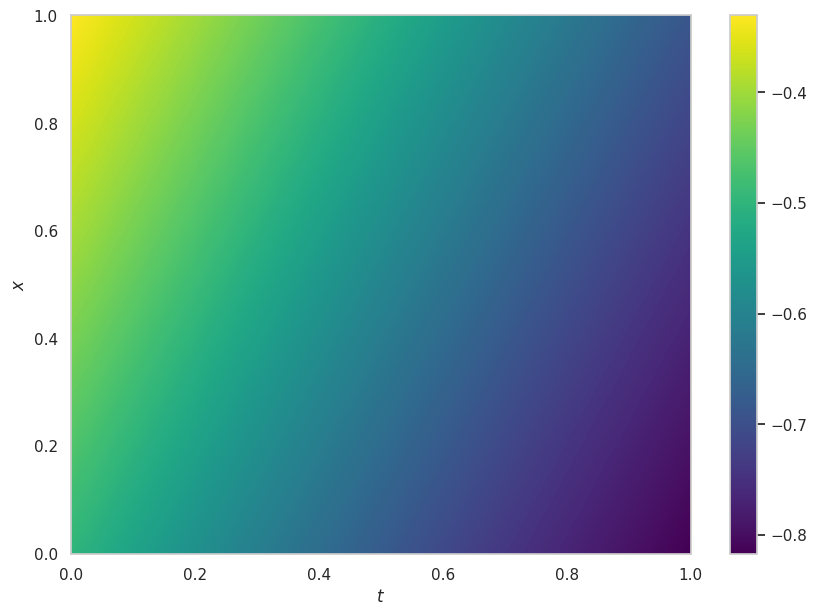} 
        \caption{PINN method}
        \label{Heat_map_of_PINN_Allen_chan.png}
    \end{subfigure}
    \caption{Heat maps plots representing the surface solutions for Problem 2 using both the exact method (left) and the PINN method (right)}
    \label{Heat map_Allen_cahn}
    \end{figure}

\begin{figure}[htbp]
    \centering
    \begin{subfigure}[b]{0.45\textwidth}
        \centering
        \includegraphics[width=\textwidth]{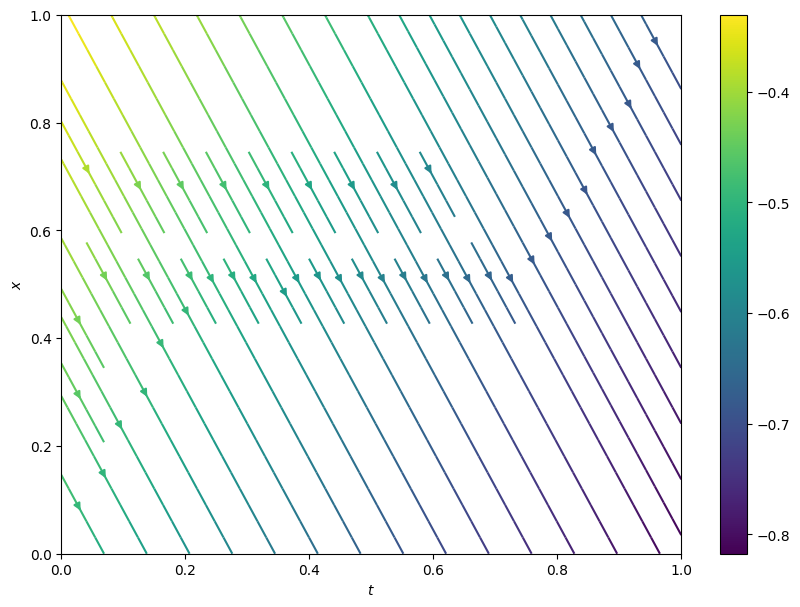} 
        \caption{Exact method}
        \label{Stream lines of exact_Allan chan.png}
    \end{subfigure}
    \hfill
    \begin{subfigure}[b]{0.45\textwidth}
        \centering
        \includegraphics[width=\textwidth]{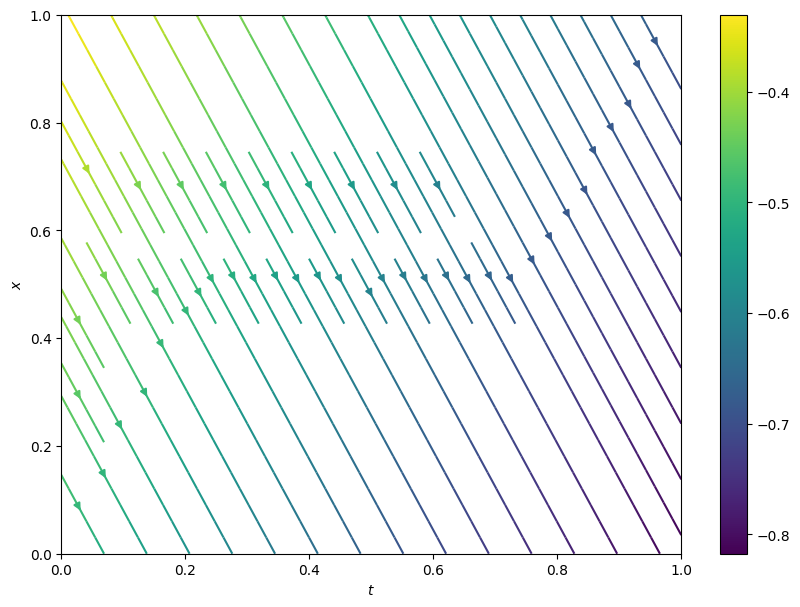} 
        \caption{PINN method}
        \label{Stream lines of PINN Allen chan.png}
    \end{subfigure}
    \caption{Gradient lines plots representing the flow of gradient vectors at solution surface. where gradient lines of exact solution is on left side and the PINN method is on right side }
    \label{streamline map_ Allen cahn}
    \end{figure}

\begin{figure}[htbp]
    \centering
    \begin{subfigure}[b]{0.45\textwidth}
        \centering
        \includegraphics[width=\textwidth]{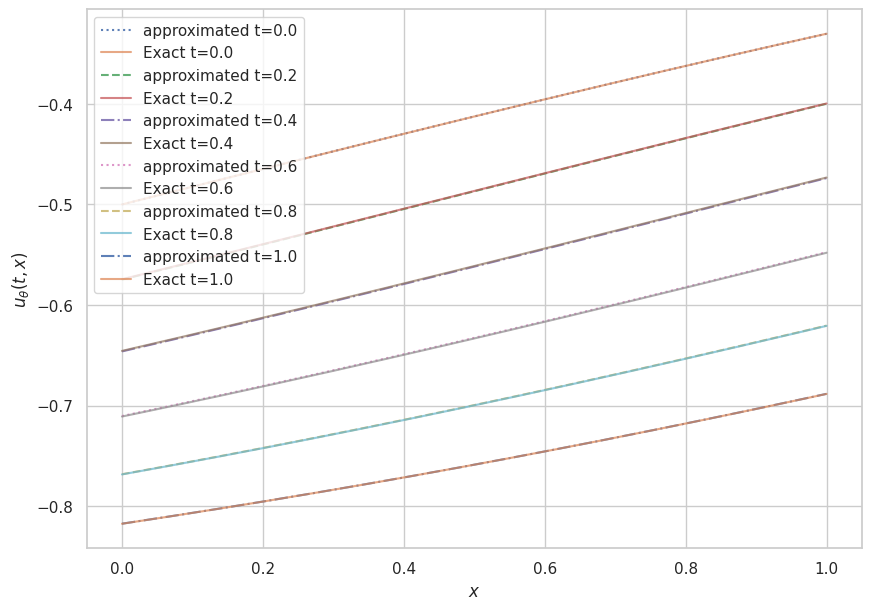} 
        \caption{Cross-section curves}
        \label{Cross_section curves Allen chan.png}
    \end{subfigure}
    \hfill
    \begin{subfigure}[b]{0.45\textwidth}
        \centering
        \includegraphics[width=\textwidth]{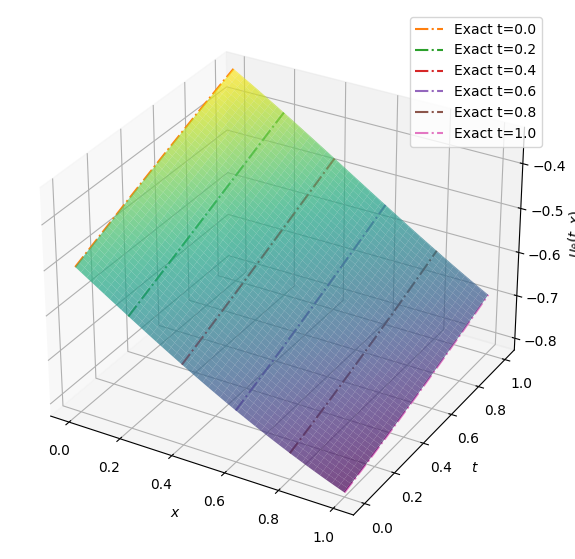} 
        \caption{Absolute error surface}
        \label{L2 norm loss_Allen_Cahn}
    \end{subfigure}
    \caption{ Comparison of cross-section curves (problem 2) between the exact method and the PINN method for various values of $t$ = 0, 0.2, 0.6, 0.8, 1. Dotted lines represent the exact solution cross-section curves, which fit well on the 3D solution surface generated by the PINN method (right). The 2D visualization of it, is shown on the left.}
    \label{loss_Allen_Cahn}
    \end{figure}

\begin{figure}[htbp]
    \centering
    \begin{subfigure}[b]{0.45\textwidth}
        \centering
        \includegraphics[width=\textwidth]{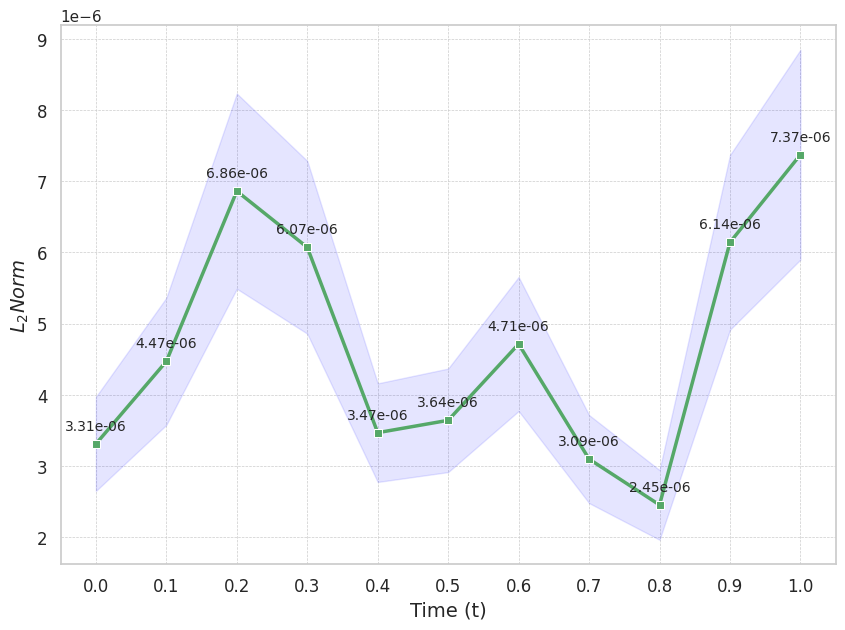} 
        \caption{\(L_2\) norm}
        \label{L_2 norm Allen cahn}
    \end{subfigure}
    \hfill
    \begin{subfigure}[b]{0.45\textwidth}
        \centering
        \includegraphics[width=\textwidth]{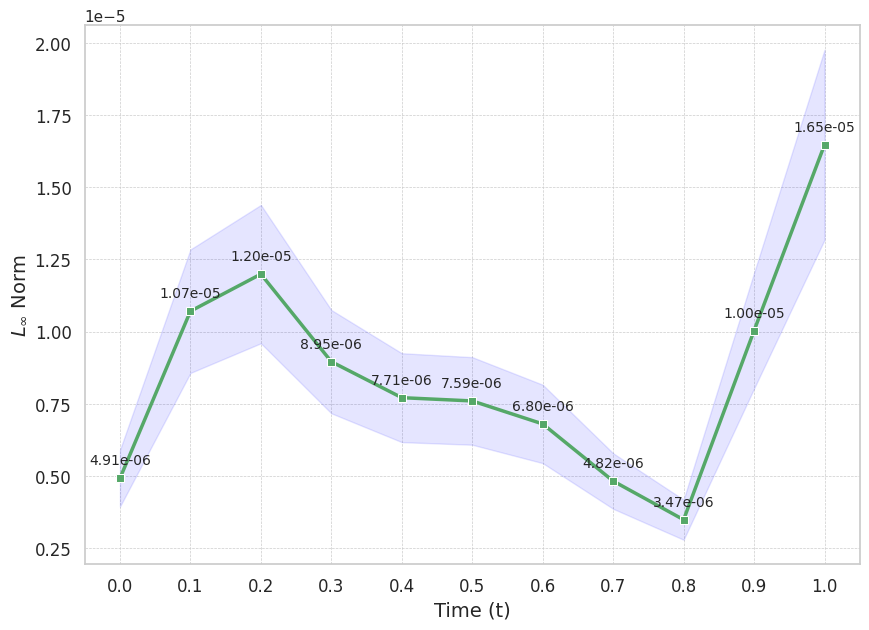} 
        \caption{$L_{\infty} norm$}
        \label{L_infty norm Allen cahn}
    \end{subfigure}
    \caption{ It represents the \(L_2\) norm and $L_{\infty}$ norm for problem 2 at t = 0, 0.1, 0.2, 0.3, 0.4, 0.5, 0.6, 0.7, 0.8, 0.9, 1  by using PINN, where $  h = 0.001$}.
    \label{L2 and L_infinity Allen cahn}
    \end{figure}

\newpage
\section{Summary and conclusion}

This research focuses on solving the fundamental models of PDEs, specifically the NWS equation and the Allen-Cahn model, using PINNs, and comparing the results with Spline numerical techniques. The study demonstrates that both models are solved very efficiently using the PINN approach. For the NWS equation, PINN consistently outperforms three variants of the Splines technique in terms of accuracy at all points mentioned in the literature. Similarly, for the Allen-Cahn model, PINN performs better than the Spline variants at all points, except at the boundary points.
\\
Additionally, the computation time for evaluating the solutions was measured, ranging from 1,000 to 10,000 points with intervals of 1,000 points. The results show that the model exhibits linear behavior in terms of computational complexity, i.e.,
O(n). Specifically, the solutions for 10,000 domain points can be obtained in approximately 107 seconds for the NWS equation and 104 seconds for the Allen-Cahn model using the PINN approach. These findings highlight the computational efficiency and accuracy of PINNs compared to traditional Spline numerical techniques.

\bibliographystyle{unsrt} 
\bibliography{main} 
\end{document}